\title{The family of ternary cyclotomic polynomials with one free prime}
\author{Yves Gallot, Pieter Moree and Robert Wilms}
\documentclass[12pt]{article}
\usepackage{amssymb, latexsym, amsfonts,amsmath,paralist}
\def\@ptsize{2}
\newtheorem{Thm}{Theorem}
\newtheorem{Con}{Conjecture}
\newtheorem{Prob}{Problem}
\newtheorem{Lem}{Lemma}

\newtheorem{Qu}{Question}

\newtheorem{cor}{Corollary}

\newtheorem{Prop}{Proposition}
\newcommand{\qed}{\hfill $\Box$}

\begin{document}
\date{}
\maketitle
{\def\thefootnote{}
\begin{abstract}
\noindent A cyclotomic polynomial $\Phi_n(x)$ is said to be ternary if $n=pqr$ with $p,q$ and $r$
distinct odd primes. Ternary cyclotomic polynomials are the simplest ones for which the
behaviour of the coefficients is not completely understood. Here we establish some results and
formulate some conjectures regarding the coefficients
appearing in the polynomial family $\Phi_{pqr}(x)$ with $p<q<r$, $p$ and $q$ fixed and $r$ a free prime.
\end{abstract}
\section{Introduction}
The $n$-th cyclotomic polynomial $\Phi_n(x)$ is defined by
$$\Phi_n(x)=\prod_{1\le j\le n\atop (j,n)=1}(x-\zeta_n^j)=\sum_{k=0}^{\infty}a_n(k)x^k,$$ with
$\zeta_n$ a $n$-th primitive root of unity (one can take $\zeta_n=e^{2\pi i/n}$).
It has degree $\varphi(n)$, with $\varphi$ Euler's totient function. 
We write 
$A(n)=\max\{|a_n(k)|:k\ge 0\}$, and this quantity is called the height of $\Phi_n(x)$. It is easy
to see that $A(n)=A(N)$, with $N=\prod_{p|n,~p>2}p$ the odd squarefree kernel. In deriving this,
 one uses the observation that if $n$ is odd, then $A(2n)=A(n)$. If $n$ has at most two distinct odd prime factors, 
then $A(n)=1$. If $A(n)>1$, then we necessarily must have
that $n$ has at least three distinct odd prime factors. In particular for $n<105=3\cdot 5\cdot 7$ we have $A(n)=1$. It turns
out that $A(105)=2$ with $a_{105}(7)=-2$. Thus the easiest case where we can expect non-trivial
behaviour of the coefficients of $\Phi_n(x)$ is the ternary case, where $n=pqr$, with $2<p<q<r$ odd primes.
In this paper we are concerned with the family of ternary cyclotomic polynomials
\begin{equation}
\label{family}
\{\Phi_{pqr}(x)|r>q\},
\end{equation}
where $2<p<q$ are fixed primes and $r$ is a `free prime'. Up to now in the literature the above family
was considered, but with also $q$ free. The maximum coefficient (in absolute value) that occurs in that
family will be denoted by $M(p)$, thus $M(p)=\max\{A(pqr):p<q<r\}$, with $p>2$ fixed. Similarly we define
$M(p;q)$ to be the maximum coefficient (in absolute value) that occurs in the
family (\ref{family}), thus $M(p;q)=\max\{A(pqr):r>q\}$, with $2<p<q$ fixed primes.\\

\noindent {\tt Example}. Bang \cite{Bang} proved that $M(p)\le p-1$. Since $a_{3\cdot 5\cdot 7}(7)=-2$ we infer
that $M(3)=2$. Using $a_{105}(7)=-2$ and $M(3)=2$, we infer that $M(3;5)=2$.\\

\noindent Let ${\cal A}(p;q)=\{a_{pqr}(k)|r>q,~k\ge 0\}$ be the set of coefficients occurring in the 
polynomial family (\ref{family}).
\begin{Prop}
\label{range}
We have ${\cal A}(p;q)=[-M(p;q),M(p;q)]\cap \mathbb Z$.
\end{Prop}
This shows the relevance of understanding $M(p;q)$. Let us first recall some known results concerning
the related function $M(p)$. Here we know thanks to Bachman \cite{B1}, who very slightly improved on
an earlier result in \cite{Beiter-2},
that
$M(p)\le 3p/4$. In 1968 it was conjectured by Sister Marion Beiter \cite{Beiter-1} (see 
also \cite{Beiter-2}) that
$M(p)\le (p+1)/2$. 
She proved it for $p\le 5$. Since M\"oller \cite{HM}
proved that $M(p)\ge (p+1)/2$ for $p>2$, her conjecture actually
would imply that $M(p)=(p+1)/2$ for $p>2$.
The first to show that Beiter's conjecture is false seems to have
been Eli Leher (in his PhD thesis), who gave the counter-example $a_{17\cdot29\cdot 41}(4801)=-10$,
showing that $M(17)\ge 10>9=(17+1)/2$. Gallot and Moree \cite{GM} provided 
for each $p\ge 11$ infinitely many infinitely many 
counter-examples $p\cdot q_j\cdot r_j$ with $q_j$ strictly increasing with $j$.
Moreover, they have shown that
for every $\epsilon>0$ and $p$ sufficiently large $M(p)>({2\over 3}-\epsilon)p$. They also proposed
the Corrected Beiter Conjecture: $M(p)\le 2p/3$.  The implications of their work for $M(p;q)$ are
described in Section \ref{yvespieter}.\\
\indent Proposition \ref{range} together with M\"oller's result quoted above gives a different proof of the
result, due to Bachman \cite{B2}, that $\{a_{pqr}(k)|~p<q<r\}=\Bbb Z$. For references and further results in this
direction (begun by I. Schur) see Fintzen \cite{jessica}.\\
\indent Jia Zhao and Xianke Zhang \cite{ZZ1} showed that $M(7)=4$, thus establishing the Beiter Conjecture
for $p=7$. In a later paper they eastablished the Corrected Beiter Conjecture:
\begin{Thm} {\rm Zhao and Zhang \cite{ZZ2}}.
\label{ZZ}
We have $M(p)\le 2p/3$.
\end{Thm}
This result together with some computer computation allows one to extend the list of exactly
known values of $M(p)$ (see Table 1). For a given prime $p$ by `smallest $n$', we mean
the smallest integer $n$ satisfying $A(n)=M(p)$ and with $p$ as its smallest prime divisor.  

\centerline{\bf TABLE 1}
\begin{center}
\begin{tabular}{|c|c|c|}
\hline
$p$ & $M(p)$ & smallest $n$\\
\hline
3 & 2 & $3\cdot 5\cdot 7$\\
\hline
5 & 3 & $5\cdot 7\cdot 11$\\
\hline
7 & 4 & $7\cdot 17\cdot 23$\\
\hline
11 & 7 & $11\cdot 19\cdot 601$\\
\hline
13 & 8 & $13\cdot 73\cdot 307$\\
\hline
19 & 12 & $19\cdot 53\cdot 859$\\
\hline
\end{tabular}
\end{center}
It is not known whether there is a finite procedure to determine $M(p)$. On the other hand, it is not difficult
to see that there is such a procedure for $M(p;q)$.
\begin{Prop}
\label{finitep}
Given primes $2<p<q$, there is a finite procedure to determine $M(p;q)$.
\end{Prop}

\noindent 
Recall that a set $S$ of primes is said to have {\it natural density} $\delta$ if the ratio 
$$\lim_{x\rightarrow \infty}{|\{p\le x: p\in S\}|\over \pi(x)}=\delta,$$ with $\pi(x)$ the number
of primes $p\le x$.
A further question that arises is how often the maximum value $M(p)$ is assumed. Here we have the
following theorem.
\begin{Thm}
\label{startie}
Given primes $2<p<q$, there exists a prime $q_0$ with $q_0\equiv q ({\rm mod~}p)$ and an integer $d$ such that $M(p,q)\le M(p, q_0)=M(p,q')$ for every
prime $q'\ge q_0$ satisfying $q'\equiv q_0 ({\rm mod~}d\cdot p)$. In particular the set
of primes $q$ with $M(p;q)=M(p)$ has a subset having a positive natural density.
\end{Thm}
A weaker result in this direction, namely that for a fixed prime $p\ge 11$, the set of primes $q$ 
such that $M(p;q)>(p+1)/2$ has a subset of positive natural density, follows from the work of
Gallot and Moree \cite{GM} (recall that $M(p) > (p+1)/2$ for $p\ge 11$).\\
\indent Unfortunately, the proof of Theorem \ref{startie} gives a lower bound for the density that seems to be far removed from the true 
value. In this paper
we present some constructions that allow one to obtain much better bounds for the density for small $p$. These results
are subsumed in the following main result of the paper.
\begin{Thm}
\label{exxpli}
Let $2<p\le 19$ be a prime with $p\ne 17$. Then the set of primes $q$ such that $M(p;q)=M(p)$ has a subset having natural density 
$\delta(p)$ as
given in the table below.\\
{\rm \centerline{\bf TABLE 2} \begin{center}
\begin{tabular}{|c|c|c|c|c|c|c|}
\hline
$p$ & 3 & 5 & 7 & 11 & 13 & 19 \\
\hline
$\delta(p)$ & 1 & 1 & 1 & ${2/5}$ & ${1/12}$ & ${1/9}$\\
\hline
\end{tabular}
\end{center}}
\end{Thm}
Numerical experimentation suggests that the set of primes $q$ such that $M(p;q)=M(p)$ has a natural density $\delta(p)$ as given
in the above table, except when $p=13$ in which case numerical experimentation suggests $\delta(13)=1/3$.\\
\indent In order to prove Theorem \ref{exxpli}, we will use the following theorem dealing with $2<p\le 7$. 
\begin{Thm}
\label{upto7}
For $2<p\le 7$ and $q>p$ we have $M(p;q)=(p+1)/2$, except in the case $p=7$, $q=13$ where $M(7;13)=3$.
\end{Thm}
The fact that $M(7;13)=3$ can be explained. Indeed, it turns out that if $ap+bq=1$ for small (in absolute
value) integers $a$ and $b$, then $M(p;q)$ is small. For example, one has the following result.
\begin{Thm}
\label{tweepeemineen}
If $p\ge 5$ and $2p-1$ is a prime, then $M(p;2p-1)=3$.
\end{Thm}
This result and similar ones are established in Section \ref{close}.\\
\indent Our main conjecture on $M(p;q)$ is the following one.
\begin{Con}
\label{startie1}
Given a prime $p$, there exists an integer $d$ and a function $g:(\mathbb Z/d\mathbb Z)^*\rightarrow \mathbb Z_{>0}$
such that for some $q_0>d$ we have for every prime $q\ge q_0$ that $M(p;q)=g({\overline q})$, where $1\le {\overline q}<d$ satisfies
$q\equiv {\overline q}({\rm mod~}d)$. The function $g$ is symmetric, that is we have $g(\alpha)=g(d-\alpha)$.
\end{Con}
The smallest integer $d$ with the above properties, if it exists, we call the {\it ternary conductor} ${\mathfrak f}_p$. 
The corresponding smallest choice of $q_0$ (obtained on setting $d={\mathfrak f}_p$) we call the {\it ternary minimal prime}.
For $p=7$ we obtain, e.g., ${\mathfrak f}_7=1$ and $q_0=17$ (by Theorem \ref{upto7}). Note that once
we know $q_0$ it is a finite computation to determine $d$ and the function $g$.
Theorem \ref{upto7} can be used to obtain 
the $p\le 7$ part of the following observation concerning the ternary conductor.
\begin{Prop}
\label{7up}
If $2<p\le 7$, then the ternary conductor exists and we have ${\mathfrak f}_p=1$. If $p\ge 11$ and ${\mathfrak f}_p$ exists,
then $p|{\mathfrak f}_p$.
\end{Prop}

While Theorem \ref{startie} only says that the set of primes $q$ with $M(p;q)=M(p)$ has a subset having a positive natural density, 
Conjecture \ref{startie1} implies that the set actually has a natural density in $\mathbb Q_{>0}$ which can 
be easily explicitly computed assuming we know $q_0$.
In order to establish this implication one can invoke
a quantitative form of Dirichlet's prime number theorem to the effect that, for $(a,d)=1$, we have, as $x$ tends to infinity,
\begin{equation}
\label{Dirichleti}
\sum_{p\le x,~p\equiv a({\rm mod~}d)}1\sim {x\over \varphi(d) \log x}.
\end{equation}
This result implies that asymptotically the primes are equidistributed over the primitive congruence classes
modulo $d$. (Recall that Dirichlet's prime number theorem, Dirichlet's theorem for
short, says that each primitive residue class contains infinitely many
primes.)

The main tool in this paper is Kaplan's lemma and is presented in Section \ref{nkaplan}. The material
in that section (except for Lemma \ref{gelijk} which is new), is taken from \cite{GM2}. As a demonstration
of working with Kaplan's lemma two examples (with and without table) are given in Section \ref{sample}.
In \cite{GMW}, the full version of this paper, details of further proofs using Kaplan's lemma can
be found. In the shorter version we have merely written `Apply Kaplan's lemma'.

The above summary of results makes clear how limited presently our knowledge of $M(p;q)$ is. For the benefit
of the interested reader we present a list of open problems in the final section.

\section{Proof of two propositions and Theorem \ref{startie}}
\label{rw}
\noindent {\it Proof of Proposition} \ref{range}. By the definiton of $M(p;q)$ we have
$${\cal A}(p;q)\subseteq [-M(p;q),M(p;q)]\cap \mathbb Z.$$ 
Let $r>q$ be a prime
such that $A(pqr)=M(p;q)$ and suppose w.l.o.g. that $a_{pqr}(k)=M(p;q)$. 
Gallot and Moree \cite{GM2} showed that we have $|a_n(k)-a_n(k-1)|\le 1$ 
for ternary $n$ (see 
Bachman \cite{B4} and Bzd{\c e}ga \cite{BZ} for alternative proofs).
Since $a_{pqr}(k)=0$
for every $k$ large enough, it then follows that $0,1,\ldots,M(p;q)$ are in ${\cal A}(p;q)$.
By a result of Kaplan \cite{Kaplan} (see Zhao and Zhang \cite{ZZ1} for a different proof), we can find a prime $s\equiv -r({\rm mod~}pq)$
and an integer $k_1$ such that $a_{pqs}(k_1)=-M(p;q)$. By a similar arguments as above one then
infers that $-M(p;q),-M(p;q)+1,\ldots,-1,0$ are all in ${\cal A}(p;q)$. \qed\\

\noindent {\it Proof of Proposition} \ref{finitep}. Let ${\cal R}_{pq}$ be a set of primes, all exceeding $q$ such that
every primitive residue class modulo $pq$ is represented. By \cite[Theorem 2]{Kaplan} we have $A(pqr)=A(pqs)$ if
$s\equiv r({\rm mod~}pq)$ with $s,r$ both primes exceeding $q$ and hence
$$M(p;q)=\max\{A(pqr):r\in {\cal R}_{pq}\}.$$
Since the computation of ${\cal R}_{pq}$ and $A(pqr)$ is a finite one, the computation of $M(p;q)$ is
also finite. \qed\\

\noindent The remainder of the section is devoted to the proof of Theorem \ref{startie}.\\
\indent For coprime positive (not necessary prime) integers $p,q,r$ we define
$$\Phi'_{p,q,r}(x)=\frac{(x^{pqr}-1)(x^p-1)(x^q-1)(x^r-1)}{(x-1)(x^{pq}-1)(x^{pr}-1)(x^{qr}-1)}=
\sum_{k=0}^{\infty}a'_{p,q,r}(k) x^k.$$
Here we do not assume $p<q<r$. Hence we have the symmetry $\Phi'_{p,q,r}(x)=\Phi'_{p,r,q}(x)$.
A routine application of the the inclusion-exclusion principle to the roots of the factors shows
that $\Phi'_{p,q,r}(x)$ is a polynomial. It is referred to as a ternary inclusion-exclusion polynomial. Inclusion-exclusion
polynomials can be defined in great generality, and the reader is referred to Bachman \cite{B4}
for an introductory discussion. He shows that such polynomials and thus $\Phi'_{p,q,r}(x)$ in particular, can be
written as products of cyclotomic polynomials (\cite[Theorem 2]{B4}).\\
\indent Analogously to $A(pqr)$ and $M(p;q)$ we define the following quantities:
$$A'(p,q,r)=\max\{|a'_{p,q,r}(k)|: k\ge 0\} , M'(p;q)=\max\{A'(p,q,r): r\ge 1\} $$
$$\mbox{ and } M'(p)=\max\{M'(p;q): q\ge 1\}.$$
We have $\Phi_{pqr}(x)=\Phi'_{p,q,r}(x)$ if $p,q,r$ are distinct primes and hence $A(pqr)=A'(p,q,r)$ in this case.

\begin{Lem}
\label{Qkaplan} For coprime positive (not necessary prime) integers $p,q,r$ we have
$A'(p,q,r_1)\le A'(p,q,r_2)\le A'(p,q,r_1)+1$ if $r_2\equiv r_1 ({\rm mod~}pq)$ and $r_2>r_1$.
\end{Lem}
{\it Proof}. Note that $r_2>\max\{p,q\}$. If $r_1>\max\{p,q\}$, then Kaplan, cf. proof of Theorem 2 in \cite{Kaplan}, 
showed that $A'(p,q,r_1)=A'(p,q,r_2)$. In the remaining case $r_1<\max\{p,q\}$, we have 
$A'(p,q,r_1)\le A'(p,q,r_2)\le A'(p,q,r_1)+1$ by the Theorem in \cite{BM}. \qed\\

\noindent In Bachman and Moree \cite{BM} it is remarked that 
$A'(p,q,r_2)=A'(p,q,r_1)+1$ can occur.
\begin{Lem}
\label{MM}
If $p$ is a prime, then $M'(p)=M(p)$. If $q$ is also a prime with $q>p$ then $M'(p;q)=M(p;q)$.
\end{Lem}
\noindent {\it Proof}.
Let $p<q$ be primes. Assume $M'(p;q)=A'(p,q,r)$, where $r$ is not necessary a prime.
By Dirichlet's theorem we can find a prime $r'$ satisfying $r'\equiv r ({\rm mod~}pq)$ and
$r'>\max(q,r)$. Therefore we have by Lemma \ref{Qkaplan}:
$$M'(p;q)=A'(p,q,r)\le A'(p,q,r')=A(p,q,r')\le M(p;q).$$
Since obviously $M(p;q)\le M'(p;q)$, we have $M'(p;q)=M(p;q)$.\\
\indent Now let only $p$ be a prime. Assume $M'(p)=A'(p,q,r)$, where $q$ and $r$ are not necessary primes.
Again by Dirichlet's theorem we find a prime $q'$ with $q'\equiv q ({\rm mod~}pr)$ and
$q'>\max(p,q)$. Using Lemma \ref{Qkaplan} we have:
$$M'(p)=A'(p,q,r)\le A'(p,q',r)\le M'(p,q')=M(p,q')\le M(p).$$
Since obviously $M(p)\le M'(p)$, we have $M'(p)=M(p)$. \qed\\

\noindent {\it Proof of Theorem} \ref{startie}. We set $q_1:=q$. Let $r_i$ be a positive integer satisfying $M'(p;q_i)=A'(p,q_i,r_i)$.
Using Lemma \ref{Qkaplan} (note that $A'(p,q,r)$ is invariant under permutations of $p,q$ and $r$) we deduce:
$$M'(p;q_1)=A'(p,q_1,r_1)\le A'(p,q_2,r_1)\le A'(p,q_2,r_2)=M'(p,q_2),$$
where $q_2=q_1+pr_1$. By the same argument the sequence $q_1,q_2,q_3,\dots$ with
$q_{i+1}=q_i+pr_i$ satisfies:
$$M'(p;q_1)\le M'(p;q_2) \le M'(p;q_3) \le \dots$$
Since $M'(p;q)\le M'(p)=M(p)$ and by, e.g., Lemma \ref{kapel}, $M(p)$ is finite, there are only finitely many different values for $M'(p;q)$. Hence
there is an index $k$ such that $M'(p;q_k)=M'(p;q_{k+i})$ for all $i\ge 0$.
That means:
$$M'(p;q_k)=A'(p,q_k,r_k)= A'(p,q_{k+1},r_k)= A'(p,q_{k+1},r_{k+1})=M'(p,q_{k+1}),$$
and by induction $A'(p,q_{k+i},r_k)= A'(p,q_{k+i},r_{k+i})$. Therefore we can
assume $r_{k+i}=r_k$ for $i\ge0$. Then we have $q_{k+i}=q_k+i \cdot pr_k$.
We set $q_0:=q_k$ and $d:=r_k$. Certainly we have $q_0\equiv q ({\rm mod~}p)$. Let $q'\ge q_0$ be a prime with
$q'\equiv q_0 ({\rm mod~}d\cdot p)$. There must be an integer $m$ such that
$q'=q_{k+m}$. Since $M'(p;q)=M(p;q)$ by Lemma \ref{MM}, we have:
$$M(p;q_1)\le M(p;q_0) = M(p;q').$$
Applying this to $M(p;q_1)$ with $M(p;q_1)=M(p)$, where we have chosen $q_1$ such that $M(p;q_1)=M(p)$, we get infinitely many primes of 
the form $q_i=q_1+i\cdot pr_1$ satisfying
$M(p;q_i)=M(p)$. On invoking (\ref{Dirichleti}) with $a=q_1$ and $d=pr_1$ the proof is then completed.
\qed

\section{The bounds of Bachman and Bzd{\c e}ga}
\label{bBB}
Let $q^*$ and $r^*$, $0<q^*,r^*<p$ be the inverses of $q$ and $r$ modulo $p$ respectively.
Set $a=\min(q^*,r^*,p-q^*,p-r^*)$. Put $b=\max(\min(q^*,p-q^*),\min(r^*,p-r^*))$. In
the sequel we will use repeatedly that $b\ge a$.
Bachman in 2003 \cite{B1} showed that
\begin{equation}
\label{babound}
A(pqr)\le \min\left({p-1\over 2}+a,p-b\right).
\end{equation}
This was more recently improved by Bzd{\c e}ga \cite{BZ} who showed that
\begin{equation}
\label{barbound}
A(pqr)\le \min(2a+b,p-b).
\end{equation}
It is not difficult to show that $\min(2a+b,p-b)\le \min({p-1\over 2}+a,p-b)$ and thus
Bzd{\c e}ga's bound is never worse than Bachman's and in practice often strict inequality holds.\\
\indent Note that if $q\equiv \pm 1({\rm mod~}p)$, then (\ref{babound}) implies that $A(pqr)\le (p+1)/2$, a result
due to Sister Beiter \cite{Beiter-1} and, independently, Bloom \cite{Bloom}.\\
\indent We like to remark that Bachman and Bzd{\c e}ga define $b$ as follows:
$$b=\min(b_1,p-b_1),~ab_1qr\equiv 1({\rm mod~}p),~0<b_1<p.$$ It is an easy exercise to
see that our definition is equivalent with this one.\\
\indent We will show that both (\ref{babound}) and (\ref{barbound}) give rise to the same upper bound
$f(q^*)$ for $M(p;q)$. Write $q^*\equiv j({\rm mod~}p)$, $r^*\equiv k({\rm mod~}p)$ with
$1\le j,k\le p-1$. Thus the right hand side of both (\ref{babound}) and (\ref{barbound}) are functions
of $j$ and $k$, which we denote by $GB(j,k)$, respectively $BB(j,k)$. We have
$$BB(j,k)=\min(2a+b,p-b)\le \min\left({p-1\over 2}+a,p-b\right)=GB(j,k),$$
with $a=\min(j,k,p-j,p-k)$ and $b=\max(\min(j,p-j),\min(k,p-k))$.
\begin{Lem}
\label{upper0}
Let $1\le j\le p-1$. Denote $GB(j,j)$ by $f(j)$. We have
$$\max_{1\le k\le p-1}BB(j,k)=\max_{1\le k\le p-1}GB(j,k)=f(j),~{\rm with}$$
$$f(j)=\begin{cases}(p-1)/2+j   & \mbox{ if } j<p/4;\\  p-j &     \mbox{ if } p/4<j\le (p-1)/2,
\end{cases}$$
and $f(p-j)=f(j)$  if $j>(p-1)/2$.
\end{Lem}
{\it Proof}. Since the problem is  symmetric under replacing $j$ by $p-j$, w.l.o.g. we may assume
that $j\le (p-1)/2$. If $j<p/4$, then
$$GB(j,k)\le {p-1\over 2}+a\le {p-1\over 2}+j=GB(j,j).$$
If $j>p/4$, then
$$GB(j,k)\le p-b\le p-j=GB(j,j).$$
Note that
$$GB(j,j)=\begin{cases}BB(j,{p+1\over 2}-j) & \mbox{ if } j<p/4; \\ BB(j,j) & \mbox{ if } j>p/4.\end{cases}$$
(E.g., if $j<p/4$, then the choice $q^*=j$, $r^*=(p+1)/2-j$ leads to $a=j$ and $b=(p+1)/2-j$ and hence
$BB(j,(p+1)/2-j)=\min((p+1)/2+j,(p-1)/2+j)=GB(j,j)$.)
Since $BB(j,k)\le GB(j,k)\le GB(j,j)$ we are done. \qed

\begin{Thm}
\label{upper1}
Let $2<p<q$. Then $M(p;q)\le f(q^*)$.
\end{Thm}
{\it Proof}. By (\ref{barbound}) and the definition of $BB(j,k)$ we have
$$M(p;q)\le \max_{1\le k\le p-1}BB(q^*,k)=f(q^*),$$ 
completing the proof. \qed\\

\noindent Lemma \ref{upper0} shows that using either (\ref{babound}) or (\ref{barbound}), we cannot improve on the upper
bound given in Theorem \ref{upper1}. Since 
$$\max_{1\le j\le p-1}f(j)=p-1-\left[{p\over 4}\right]=
\begin{cases}
3(p-1)/4 & {\rm if~} p\equiv 1({\rm mod~}4);\cr 
(3p-1)/4 & {\rm if~} p\equiv 3({\rm mod~}4),
\end{cases}
$$
we infer that
$$M(p)\le \max_{1\le j\le p-1}\max_{1\le k\le p-1}GB(j,k)=\max_{1\le j\le p-1}f(j)<{3\over 4}p.$$

\section{Earlier work on $M(p;q)$}
\label{yvespieter}
Implicit in the literature are various results on $M(p;q)$ (although we are the first to explicitly study
$M(p;q)$). Most of these are mentioned in the rest of this paper. Here we rewrite the main result of
Gallot and Moree \cite{GM} in terms of $M(p;q)$ and use it for $p=11$ and $p=13$ (to deal with
$q\equiv 4({\rm mod~}11)$, respectively $q\equiv 5({\rm mod~}13)$).
\begin{Thm}
\label{oud}
Let $p\ge 11$ be a prime. Given any $1\le \beta \le p-1$ we let $\beta^*$ be the unique integer $1\le \beta^* \le p-1$ with
$\beta \beta^* \equiv 1({\rm mod~}p)$. Let ${\cal B}_{-}(p)$ be the set of integers satisfying
$$1\le \beta \le {p-3\over 2},~p\le \beta+2\beta^*+1,~\beta>\beta^*.$$
Let ${\cal B}_{+}(p)$ be the set of integers satisfying
$$1\le \beta \le {p-3\over 2},~p\le\beta+\beta^*,~\beta\ge \beta^*/2.$$
Let ${\cal B}(p)$ be the union of these (disjoint) sets. As $(p-3)/2\in {\cal B}(p)$, it is non-empty. 
Let $q\equiv \beta ({\rm mod~}p)$ be a prime satisfying $q>p$.
Suppose that the inequality $q> q_{-}(p):=p(p-\beta^*)(p-\beta^*-2)/(2\beta)$ holds if $\beta \in {\cal B}_{-}(p)$ and
$$q> q_{+}(p):={p(p-1-\beta)\over \gamma(p-1-\beta)-p+1+2\beta},$$
with $\gamma=\min((p-\beta^*)/(p-\beta),(\beta^*-\beta)/\beta^*)$ if $\beta\in {\cal B}_{+}(p)$.
Then
$$M(p;q)\ge p-\beta >{p+1\over 2}$$
and hence $M(p)\ge p-\min\{{\cal B}(p)\}$.
\end{Thm}
We have ${\cal B}(11)=\{4\}, {\cal B}(13)=\{5\}, {\cal B}(17)=\{7\}$ and ${\cal B}(19)=\{8\}$. In general one
can show \cite{CGM} using Kloosterman sum techniques that
$$\Big||{\cal B}(p)|-{p\over 16}\Big|\le 8\sqrt{p}(\log p+2)^3.$$
The lower bound for $M(p)$ resulting from this theorem, $p-\min\{{\cal B}(p)\}$, never exceeds $2p/3$ and this together with extensive
numerical experimentation led Gallot and Moree \cite{GM} to propose the corrected Beiter conjecture, now proved by Zhao
and Zhang (Theorem \ref{ZZ}).\\
\indent Under the appropriate conditions on $p$ and $q$, Theorem \ref{oud} says that $M(p;q)\ge p-\beta$, whereas
Theorem \ref{upper1} yields $M(p;q)\le f(\beta^*)$. Thus studying the case $p-\beta=f(\beta^*)$ with $\beta\in {\cal B}(p)$, leads
to a small subset of cases where $M(p;q)$ can be exactly computed using Theorem \ref{oud}.
\begin{Thm}
Let $p\ge 13$ with $p\equiv 1({\rm mod~}4)$ be a prime. Let $x_0$ be the smallest positive integer
such that $x_0^2+1\equiv 0({\rm mod~}p)$. If $x_0>p/3$, $q\equiv x_0({\rm mod~}p)$ and $q\ge q_{+}(p)$ (with
$\beta=x_0$), then $M(p;q)=p-x_0$. 
\end{Thm}
{\it Proof}. Some easy computations show that if $p-\beta=f(\beta^*)$ and $\beta\in {\cal B}(p)$, we must 
have $\beta\in {\cal B}_{+}(p)$, ${p-1\over 2}<\beta^*<{3\over 4}p$ and hence $f(\beta^*)=\beta^*$ and so
\begin{equation}
\label{beeta}
~\beta\in {\cal B}_{+}(p),~1\le \beta \le {p-3\over 2},~\beta+\beta^*=p,~\beta^*\le 2\beta,~{p-1\over 2}<\beta^*<{3\over 4}p.
\end{equation}
Note that $\beta+\beta^*=p$, $p\ge 13$, has a solution with $\beta<p/2$ iff $p\equiv 1({\rm mod~}4)$ and
$\beta=x_0$ (and hence $\beta^*=p-x_0$) with $x_0$ the smallest solution of $x_0^2+1 \equiv 0({\rm mod~}p)$. If
$x_0>p/3$, then $\beta=x_0$ satisfies (\ref{beeta}). Since by assumption 
$q\ge q_{+}(p)$ and $q\equiv x_0({\rm mod~}p)$, we have $M(p;q)\ge p-x_0$ by Theorem \ref{oud}. On the
other hand, by Theorem \ref{upper1}, we have $M(p;q)\le f(p-x_0)=f(x_0)=p-x_0$.\qed\\

\noindent {\tt Remark}. The set of primes $p$ satisfying $p\equiv 1({\rm mod~}4)$ and
$x_0>p/3$ (which starts $\{13, 29, 53, 73, 89, 173, \cdots\}$) has natural density $1/6$. This follows on taking 
$\alpha_2=1/2$ and $\alpha_1=1/3$ in the result of Duke et al. \cite{DFI}, that if 
$f$ is a quadratic polynomial with complex roots and $0\le \alpha_1 < \alpha_2\le 1$ are prescribed
real numbers, then as $x$ tends to infinity,
$$\#\{(p,v):p\le x,~f(v)\equiv 0({\rm mod~}p),~\alpha_1\le {v\over p}<\alpha_2\}\sim (\alpha_2-\alpha_1)\pi(x).$$

\section{Computation of $M(3;q)$}
Note that for all primes $q$ and $r$ with $1<q<r$, there exists some unique $h\le (q-1)/2$ and $k>0$ such
that
$r=(kq+1)/h$ or $r=(kq-1)/h$. If $n\equiv 0({\rm mod~}3)$ is
ternary, then either $A(n)=1$ or $A(n)=2$ as $M(3)=2$. The following result due to Sister Beiter \cite{Beiter-3}
allows one to compute $A(n)$ in this case.
\begin{Thm} 
\label{b}
Let $n\equiv 0({\rm mod~}3)$ be ternary.\\
If $h=1$, then $A(n)=1$ iff $k\equiv 0({\rm mod~}3)$.\\
If $h>1$, then $A(n)=1$ iff one of the following conditions holds:\\
{\rm (a)} $k\equiv 0({\rm mod~}3)$ and $h+q\equiv 0({\rm mod~}3)$.\\
{\rm (b)} $k\equiv 0({\rm mod~}3)$ and $h+r\equiv 0({\rm mod~}3)$.
\end{Thm}
We have seen that $M(3;5)=2$. The next result extends this.
\begin{Thm}
\label{driedrie}
Let $q>3$ be a prime. We have $M(3;q)=2$.
\end{Thm}
{\it Proof}. In case $q\equiv 1({\rm mod~}3)$, then let $r$ be a prime such that
$r\equiv 1+q({\rm mod~}3q)$. Since $(1+q,3q)=1$, there are in fact infinitely many
such primes (by Dirichlet's theorem). In case $q\equiv 2({\rm mod~}3)$, then let $r$ be a prime such that
$r\equiv 1+2q({\rm mod~}3q)$. Since $(1+2q,3q)=1$, there are infinitely many
such primes. The prime $r$ was chosen so to ensure that $h=1$ and $3\nmid k$. Using Theorem
\ref{b} it then follows that $A(3qr)=2$ and hence $M(3;q)=2$. \qed
\section{Kaplan's lemma reconsidered}
\label{nkaplan}
\indent Our main tool will be the following recent result due to Kaplan \cite{Kaplan}, the
proof of which uses the identity
$$\Phi_{pqr}(x)=(1+x^{pq}+x^{2pq}+\cdots)(1+x+\cdots+x^{p-1}-x^q-\cdots-x^{q+p-1})
\Phi_{pq}(x^r).$$
\begin{Lem} {\rm (Nathan Kaplan, 2007)}.
\label{kapel}
Let $2<p<q<r$ be primes and $k\ge 0$ be an integer.
Put $$b_i=\begin{cases}a_{pq}(i) & \mbox{ if }ri\le k; \\ 0 & \mbox{ otherwise.}\end{cases}$$ We have
\begin{equation}
\label{lacheens}
a_{pqr}(k)=\sum_{m=0}^{p-1}(b_{f(m)}-b_{f(m+q)}),
\end{equation}
where $f(m)$ is the unique integer  such that
$f(m)\equiv r^{-1}(k-m) ({\rm mod~}pq)$ and $0\le  f(m) < pq$.
\end{Lem}
(If we need to stress the $k$-dependence of $f(m)$, we will write $f_k(m)$ instead of $f(m)$,
see, e.g., Lemma \ref{gelijk} and its proof.)
This lemma reduces the computation of $a_{pqr}(k)$ to that of $a_{pq}(i)$ for
various $i$. These binary cyclotomic polynomial coefficients are computed
in the following lemma. For a proof
see, e.g., Lam and Leung \cite{LL} or Thangadurai \cite{Thanga}.
\begin{Lem}
\label{binary}
Let $p<q$ be odd primes. Let $\rho$ and $\sigma$ be the (unique) non-negative
integers for which $1+pq=(\rho+1) p+(\sigma+1) q$.
Let $0\le m<pq$. Then either $m=\alpha_1p+\beta_1q$ or $m=\alpha_1p+\beta_1q-pq$
with $0\le \alpha_1\le q-1$ the unique integer such that $\alpha_1 p\equiv m({\rm mod~}q)$
and $0\le \beta_1\le p-1$ the unique integer such that $\beta_1 q\equiv m({\rm mod~}p)$.
The cyclotomic coefficient $a_{pq}(m)$ equals
$$\begin{cases}1 & \mbox{ if }m=\alpha_1p+\beta_1q \mbox{ with }0\le \alpha_1\le \rho,~0\le \beta_1\le
\sigma;\\ -1 & \mbox{ if }m=\alpha_1p+\beta_1q-pq \mbox{ with }\rho+1\le \alpha_1\le q-1,~\sigma+1\le 
\beta_1\le p-1;\\  0 & \mbox{ otherwise.}\end{cases}$$
\end{Lem}
We say that $[m]_p=\alpha_1$ is the {\it $p$-part of $m$} and $[m]_q=\beta_1$ is the {\it $q$-part
of $m$}. It is easy to see that
$$m=\begin{cases}[m]_pp+[m]_qq & \mbox{ if }[m]_p\le \rho \mbox{ and }[m]_q\le \sigma;\\ [m]_pp+[m]_qq-pq & \mbox{ if }[m]_p>\rho \mbox{ and }[m]_q>\sigma;\\ [m]_pp+[m]_qq-\delta_mpq & \mbox{ otherwise,}\end{cases}$$
with $\delta_m\in \{0,1\}$. Using this observation we find that, for $i<pq$,
$$b_i=\begin{cases}1 & \mbox{ if }[i]_p\le \rho, [i]_q\le \sigma \mbox{ and }[i]_pp+[i]_qq\le k/r;\\ -1 & \mbox{ if }[i]_p>\rho, [i]_q>\sigma \mbox{ and }[i]_pp+[i]_qq-pq\le k/r;\\ 0 & \mbox{ otherwise.}\end{cases}$$
Thus in order to evaluate $a_{pqr}(n)$ using Kaplan's lemma it suffices to compute $[f(m)]_p$, $[f(m)]_q$, 
and $[f(m+q)]_q$ (note that $[f(m)]_p=[f(m+q)]_p$).\\
\indent For future reference we provide a version of Kaplan's lemma in which the
computation of $b_i$ has been made explicit, and thus is self-contained.
\begin{Lem}
\label{kapel2}
Let $2<p<q<r$ be primes and $k\ge 0$ be an integer.
We put $\rho=[(p-1)(q-1)]_p$ and $\sigma=[(p-1)(q-1)]_q$. 
Furthermore, we put $$b_i=\begin{cases}1 & \mbox{ if }[i]_p\le \rho, [i]_q\le \sigma \mbox{ and }[i]_pp+[i]_qq\le k/r;\\ -1 & \mbox{ if }[i]_p>\rho, [i]_q>\sigma \mbox{ and }[i]_pp+[i]_qq-pq\le k/r;\\ 0 & \mbox{ otherwise.}\end{cases}$$
We have
\begin{equation}
\label{lacheens2}
a_{pqr}(k)=\sum_{m=0}^{p-1}(b_{f(m)}-b_{f(m+q)}),
\end{equation}
where $f(m)$ is the unique integer  such that
$f(m)\equiv r^{-1}(k-m) ({\rm mod~}pq)$ and $0\le  f(m) < pq$.
\end{Lem}
Note that if $i$ and $j$ have the same $p$-part, then $b_ib_j\ne -1$, that is $b_i$ and $b_j$
cannot be of opposite sign. {}From this it follows that $|b_{f(m)}-b_{f(m+q)}|\le 1$, and
thus we infer from Kaplan's lemma that $|a_{pqr}(k)|\le p$ and hence $M(p)\le p$.\\
\indent Using the mutual coprimality of $p,q$ and $r$ we arrive at the following trivial, but useful, lemma.
\begin{Lem} \label{flow} We have $\{[f(m)]_q:0\le m\le p-1\}=\{0,1,2,\ldots,p-1\}$ and
$|\{[f(m)]_p:0\le m\le p-1\}|=p$. The same conclusions hold if we replace
$[f(m)]_q$ and $[f(m)]_p$ by $[f(m+q)]_q$, respectively $[f(m+q)]_p$.
\end{Lem}

\indent On working with Kaplan's lemma one first computes $a_{pq}(f(m))$ and then
$b_{f(m)}$. As a check on the correctness of the computations we note that the following
identity should be satisfied.
\begin{Lem}
\label{gelijk}
We have $$\sum_{m=0}^{p-1}a_{pq}(f_k(m))=\sum_{m=0}^{p-1}a_{pq}(f_k(m+q)).$$
\end{Lem}
{\it Proof}. Choose an integer $k_1\equiv k({\rm mod~}pq)$ such that $k_1>pqr$. Then
$a_{pqr}(k_1)=0$. By Lemma \ref{kapel} we find that
$$0=a_{pqr}(k_1)=\sum_{m=0}^{p-1}[a_{pq}(f_{k_1}(m))-a_{pq}(f_{k_1}(m+q))].$$
Since $f_k(m)$ only depends on the congruence class of $k$ modulo $pq$, $f_{k_1}(m)=f_k(m)$
and the result follows. \qed
\subsection{Working with Kaplan's lemma: examples}
\label{sample}
In this section we carry out some sample computations using Kaplan's lemma. 
For more involved examples the
reader is referred to \cite{GM}.\\
\indent We remark that the result that $a_n(k)=(p+1)/2$ in Lemma \ref{moellerext} is due
to Herbert M\"oller \cite{HM}. The proof we give here of this is rather different. The foundation for M\"oller's result
is due to Emma Lehmer \cite{Emma}, who already in 1936 had shown that $a_n({1\over 2}(p-3)(qr+1))=(p-1)/2$ with
$p,q,r$ and $n$ satisfying the conditions of Lemma \ref{moellerext}.
\begin{Lem}
\label{moellerext}
Let $p<q<r$ be primes satisfying
$$p>3,~q\equiv 2({\rm mod~}p),~r\equiv {p-1\over 2}({\rm mod~}p),~r\equiv {q-1\over 2}({\rm mod~}q).$$
For $k=(p-1)(qr+1)/2$ we have $a_{pqr}(k)=(p+1)/2$.
\end{Lem}
{\it Proof} (taken from \cite{GM2}.) Using that $q\equiv 2({\rm mod~}p)$, we infer from
$1+pq=(\rho+1)p+(\sigma+1)q$ that $\sigma={p-1\over 2}$ and $(\rho+1)p=1+({p-1\over 2})q$ (and
hence $\rho=(p-1)(q-2)/(2p)$). On invoking the Chinese remainder theorem one checks that
\begin{equation}
\label{LD}
-{r^{-1}}\equiv 2\equiv -\left({q-2\over p}\right)p+q({\rm mod~}pq).
\end{equation}
Furthermore, writing $f(0)$ as a linear combination of $p$ and $q$ we see that
\begin{equation}
\label{efnul}
f(0)\equiv {k\over r}\equiv \left({p-1\over 2}\right)q+{p-1\over 2r}\equiv \left({p-1\over 2}\right)q+1-p\equiv \rho p({\rm
mod~}pq).
\end{equation}
Since $f(m)\equiv f(0)-{m\over r}({\rm mod~}pq)$ we find using (\ref{LD}), (\ref{efnul}) 
and the observation that $\rho-m(q-2)/p\ge 0$ for $0\le m\le (p-1)/2$, that $[f(m)]_p=\rho-m(q-2)/p\le \rho$ and $[f(m)]_q=m\le \sigma$ 
for $0\le m\le (p-1)/2$. Since
$[f(m)]_pp+[f(m)]_qq=\rho p+2m\le \rho p+ p-1=[k/r]$, we deduce that
$a_{pq}(f(m))=b_{f(m)}=1$ in this range (see also Table 3).\\
\medskip
\begin{center}
{\bf TABLE 3}
\begin{tabular}{|c|c|c|c|c|c|c|}
\hline
$m$ & $[f(m)]_p$ & $[f(m)]_q$ & $f(m)$ & $a_{pq}(f(m))$ & $b_{f(m)}$\\
\hline
0 & $\rho$ & $0$ & $\rho p$ & 1 & 1\\
\hline
1 & $\rho-(q-2)/p$ & 1 & $\rho p+2$ & 1 & 1\\
\hline
$\vdots$ & $\vdots$ & $\vdots$ & $\vdots$ & 1 & 1\\
\hline
$j$ & $\rho-j(q-2)/p$ & $j$ & $\rho p +2j$ & 1 & 1\\
\hline
$\vdots$ & $\vdots$ & $\vdots$ & $\vdots$ & 1 & 1\\
\hline
$(p-1)/2$ & $0$ & $(p-1)/2$ & $(p-1)q/2$ & 1 & 1\\
\hline
\end{tabular}
\end{center}

Note that $f(m)\equiv f(0)-m/r\equiv \rho p+2m({\rm mod~}pq)$, from which one easily
infers that $f(m)=\rho p+2m$ for $0\le m\le p-1$ (as $\rho p+2m\le \rho p+2(p-1)<pq$). In
the range ${p+1\over 2}\le m\le p-1$ we have $f(m)\ge \rho p+p+1=(p-1)q/2+2>k/r$, and hence
$b_{f(m)}=0$.\\
\indent On noting that $f(m+q)\equiv f(m)-q/r\equiv f(m)+2q\equiv \rho p+2m+2q({\rm mod~}pq)$,
one easily finds, for $0\le m\le p-1$, that $f(m+q)=\rho p+2m+2q>k/r$ and hence $b_{f(m+q)}=0$.\\
\indent On invoking Kaplan's lemma one finds
$$a_{pqr}(k)=\sum_{m=0}^{p-1}b_{f(m)}-\sum_{m=0}^{p-1}{b_{f(m+q)}}={p+1\over 2}-0={p+1\over 2}.$$
This concludes the proof.\qed

\begin{Lem}
\label{wilms2}
Let $3<p<q<r$ be primes satisfying
$$q \equiv 1 ({\rm mod~}p), ~r^{-1}\equiv \frac{p+q}{2} ({\rm mod~} p q).$$
For $k=(p-1)qr/2 - pr + 2$ we have $a_{pqr}(k)=-\min(\frac{q-1}{p}+1, \frac{p+1}{2})$.
\end{Lem}

\noindent {\it Proof}. Let $0\le m\le p-1$. We have:
$$\rho=\frac{(p-1)(q-1)}{p}\mbox{ and } \sigma=0,$$
$$k \equiv 1 ({\rm mod~}p), ~k\equiv 0 ({\rm mod~}q), ~ k\equiv 2 ({\rm mod~}r),$$
so that we can compute:
$$[f(m)]_q \equiv q^{-1} r^{-1} (k-m) \equiv (1-m)/2 ({\rm mod~}p)$$
$$[f(m+q)]_q \equiv q^{-1} r^{-1} (k-m-q) \equiv -m/2 ({\rm mod~}p)$$
$$[f(m)]_p= [f(m+q)]_p \equiv p^{-1} r^{-1} (k-m) \equiv -m/2 ({\rm mod~}q).$$
This leads to:
$$[f(m)]_q = \begin{cases} (p+1-m)/2 & \mbox{for } m \mbox{ even} \\
(2p+1-m)/2 & \mbox{for } m \mbox{ odd and } m\neq 1\\ 0 & \mbox{for } m=1 \end{cases}$$
$$[f(m+q)]_q = \begin{cases} (p-m)/2 & \mbox{for } m \mbox{ odd} \\
(2p-m)/2 & \mbox{for } m \mbox{ even and } m \neq 0\\ 0 & \mbox{for } m=0 \end{cases}$$
$$[f(m)]_p=[f(m+q)]_p=\begin{cases} (q-m)/2 & \mbox{for } m \mbox{ odd} \\
(2q-m)/2 & \mbox{for } m \mbox{ even and } m \neq 0\\ 0 & \mbox{for } m=0. \end{cases}$$
We consider the following four cases:
\begin{itemize}
\item Case 1: $[f(m)]_p \le \rho$, $[f(m)]_q \le \sigma$.
In this case $m=1$. Therefore:
$$[f(m)]_p p + [f(m)]_q q = \frac{p(q-1)}{2}>\frac{k}{r}.$$
\item Case 2: $[f(m)]_p > \rho$, $[f(m)]_q > \sigma$.
This case only arises if $m$ is even and $m\ge 2$. Then we have:
$$[f(m)]_p p + [f(m)]_q q -pq = \frac{2q-m}{2}p + \frac{p+1-m}{2}q -pq$$
$$=\frac{q(p+1-m)-mp}{2}\le \frac{q(p-1)}{2} - p + \frac{2}{r}=\frac{k}{r}.$$
However, not all even $m\ge 2$ satisfy $[f(m)]_p > \rho$. For this it is
necessary that $\frac{2q-m}{2}>\frac{(p-1)(q-1)}{p}$. That means $\frac{m}{2}<\frac{q-1}{p}+1$ and by
$0<\frac{m}{2}\le\frac{p-1}{2}$ we have exactly $\min(\frac{q-1}{p},\frac{p-1}{2})$ different
values of $m$ in this case.
\item Case 3: $[f(m+q)]_p \le \rho$, $[f(m+q)]_q \le \sigma$. 
In this case we have $m=0$. Therefore:
$$[f(m+q)]_p p + [f(m+q)]_q q = 0 \le \frac{k}{r}.$$
\item Case 4: $[f(m+q)]_p > \rho$, $[f(m+q)]_q > \sigma$.
We must have $2|m$ and $m\ge 2$. We find:
$$[f(m+q)]_p p + [f(m+q)]_q q -pq = \frac{2q-m}{2}p + \frac{2p-m}{2}q -pq>\frac{k}{r}.$$
\end{itemize}
The above case analysis shows that (respectively),
$$\sum_{m=0\atop b_{f(m)}=1}^{p-1}1=0,~\sum_{m=0\atop b_{f(m)}=-1}^{p-1}1=
\min\left(\frac{q-1}{p},\frac{p-1}{2}\right),~\sum_{m=0\atop b_{f(m+q)}=1}^{p-1}1=1,~
\sum_{m=0\atop b_{f(m+q)}=-1}^{p-1}1=0.$$
Kaplan's lemma then yields
$$a_{pqr}(k)=\left(0-\min\left(\frac{q-1}{p},\frac{p-1}{2}\right)\right)-\left(1-0\right)=-\min\left(\frac{q-1}{p}+1,\frac{p+1}{2}\right).$$
\begin{Lem}
\label{wilms1}
Let $3<p<q<r$ be primes satisfying
$$q \equiv -2 (\mathrm{mod}~ p), ~r^{-1}\equiv p-2 ({\rm mod~}p q) \mbox{ and } q>p^2/2.$$
For $k=\frac{p+1}{2}(1+r (2-p+q))+r+q-rq$ we have 
$a_{pqr}(k)=-(p+1)/2$.
\end{Lem}
\noindent {\it Proof of Lemma } \ref{wilms1}. Apply Kaplan's lemma. \qed\\

\noindent {\tt Remark}. Numerical experimentation suggests that with this choice of $k$, a condition of the form $q>p^2c_1$, with
$c_1$ some absolute positive constant, is unavoidable.

\begin{Lem}
\label{wilms-1}
Let $3<p<q<r$ be primes satisfying
$$q \equiv -1 ({\rm mod~}p), ~r^{-1}\equiv \frac{p+q}{2} ({\rm mod~} p q) \mbox{ and } q\ge p^2-2p.$$
For $k=p (q-1) r/2 - rq + p-1$ we have $a_{pqr}(k)=-(p+1)/2$.
\end{Lem}

\noindent {\it Proof}. Apply Kaplan's lemma.\qed\\

\noindent {\it Proof of Proposition} \ref{7up}. The first assertion follows by Theorem \ref{upto7}, so
assume $p\ge 11$. We will argue by 
contradiction. So suppose that $p\nmid {\mathfrak f}_p$. Put $\beta=(p-3)/2$. By the
Chinese remainder theorem and Dirichlet's theorem there are infinitely many primes $q_1$ such
that $q_1\equiv 2({\rm mod~}p)$ and $q_1\equiv 1({\rm mod~}\mathfrak{f}_p)$.
Further, there are infinitely many primes $q_2$ such
that $q_2\equiv \beta({\rm mod~}p)$ and $q_2\equiv 1({\rm mod~}\mathfrak{f}_p)$.
By the definition of $\mathfrak{f}_p$ there exists an integer $c$ such that
$M(p;q)=c$ for all $q\equiv 1({\rm mod~}\mathfrak{f}_p)$ that are large enough.
However, by Lemma \ref{moellerext} we have $M(p;q_1)=(p+1)/2$ and by Theorem \ref{oud} (note that 
$\beta\in {\cal B}(p)$) we 
have $M(p;q_2)>(p+1)/2$ for all $q_2$ large enough. This contradiction shows that $p\nmid {\mathfrak f}_p$.\qed\\

\noindent The results from this section together with those from Section \ref{bBB} allow one to establish
the following theorem. In Section \ref{close} we will discuss the sharpness of the lower bounds for $q$.
\begin{Thm} 
\label{zeven}
Let $2<p<q$ be primes.\\
{\rm (a)} If $q\equiv 2({\rm mod~}p)$, then $M(p;q)=(p+1)/2$.\\
{\rm (b)} If $q\equiv -2({\rm mod~}p)$ and $q>p^2/2$, then $M(p;q)=(p+1)/2$.\\
{\rm (c)} If $q\equiv 1({\rm mod~}p)$ and $q\ge (p-1)p/2+1$, then $M(p;q)=(p+1)/2$.\\
{\rm (d)} If $q\equiv -1({\rm mod~}p)$ and $q\ge p^2-2p$, then $M(p;q)=(p+1)/2$.
\end{Thm}
{\it Proof}. By Theorem \ref{driedrie} we have $M(3;q)=2=(3+1)/2$, so assume $p>3$.\\
(a) We have $M(p;q)\ge (p+1)/2$ by Lemma \ref{moellerext}, and $M(p;q)\le f(2^*)=f((p+1)/2)=(p+1)/2$
by Theorem \ref{upper1}. \\
(b)+(c)+(d) Similar to that of part (a). Note 
that $f((-2)^*)=f((p-1)/2)=(p+1)/2$ and $f(1)=f(p-1)=(p+1)/2$. \qed\\

\noindent Using Theorem \ref{zeven} it is easy to establish the following result.
\begin{Thm} 
\label{5}
Let $q>5$ be a prime. Then $M(5;q)=3$. \end{Thm}
{\it Proof}. The proof is most compactly given by Table 4.\\ 
\medskip
\centerline{\bf TABLE 4}
\begin{center}
\begin{tabular}{|c|c|c|c|}
\hline
${\overline q}$ & $q_0$ & $M(5;q)$ & result\\
\hline
1 &  11 &  3 & Theorem \ref{zeven} (c) \\
\hline
2 & 7 & 3 & Theorem \ref{zeven} (a)\\
\hline
3 & 13 & 3 & Theorem \ref{zeven} (b)\\
\hline
4 & 19  & 3 & Theorem \ref{zeven} (d)\\
\hline
\end{tabular}
\end{center}
The table should be read as follows. {}From, e.g., the third row we read that for $q\equiv 3({\rm mod~}5)$, 
$q\ge 13$, we have that $M(5;q)=3$ by Theorem \ref{zeven} (b).\qed

\section{Computation of $M(7;q)$}
Theorem \ref{zeven} in addition with the following two lemmas allows one to compute $M(7;q)$. These lemmas concern the
computation of $M(p;q)$ with $q\equiv (p\pm 1)/2({\rm mod~}p)$.
\begin{Lem}
\label{p-1}
Let $p\ge 5$ be a prime. Let $q\ge\max(3p,p(p+1)/4)$ be a prime satisfying
$q\equiv {p-1\over 2}({\rm mod~}p)$. Let $r>q$ be a prime satisfying
$$r^{-1}\equiv {p+1\over 2}({\rm mod~}p),~r^{-1}\equiv p({\rm mod~}q).$$
For $k=p-1+r(1+q(p-1)/2-p(p+1)/2)$ we have $a_{pqr}(k)=(p+1)/2$.
\end{Lem}

\noindent {\it Proof}. Apply Kaplan's lemma. \qed

\begin{Lem}
\label{p+1}
Let $p\ge 5$ be a prime. Let $q\ge \max(3p,p(p-1)/4+1)$ be a prime satisfying
$q\equiv {p+1\over 2}({\rm mod~}p)$. Let $r>q$ be a prime satisfying
$$r^{-1}\equiv {p-1\over 2}({\rm mod~}p),~r^{-1}\equiv p({\rm mod~}q).$$
For $k=q+p-1+r(q(p-1)/2-p(p+1)/2)$ we have $a_{pqr}(k)=(p+1)/2$.
\end{Lem}

\noindent {\it Proof}. Apply Kaplan's lemma. \qed

\begin{Thm} \label{7be} $~$\\
{\rm (a)} Let $q\ge\max(3p,p(p+1)/4)$ be a prime satisfying
$q\equiv {p-1\over 2}({\rm mod~}p)$, then $(p+1)/2\le M(p;q)\le (p+3)/2$.\\
{\rm (b)} Let $q\ge \max(3p,p(p-1)/4+1)$ be a prime satisfying
$q\equiv {p+1\over 2}({\rm mod~}p)$, then $(p+1)/2\le M(p;q)\le (p+3)/2$.
\end{Thm}
{\it Proof}. Follows on noting that
$$f\Big(\big({p+1\over 2}\big)^*\Big)=f(2)={p+3\over 2}=f(p-2)=f\Big(\big({p-1\over 2}\big)^*\Big),$$
and combining Lemmas \ref{p-1} and \ref{p+1} with Theorem \ref{upper1}.\qed

\begin{Thm}
\label{7}
We have $M(7;11)=4$, $M(7;13)=3$ and for $q\ge 17$ a prime, $M(7;q)=4$.
\end{Thm}
{\it Proof}. The proof is most compactly given by a table (Table 5). Recall
that Zhao and Zhang \cite{ZZ1} proved that $M(7)\le 4$.\\

\centerline{\bf TABLE 5}
\begin{center}
\begin{tabular}{|c|c|c|c|}
\hline
${\overline q}$ & $q_0$ & $M(7;q)$ & result\\
\hline
1 &  29 &  4 & Theorem \ref{zeven} (c)\\
\hline
2 & 23 & 4 & Theorem \ref{zeven} (a)\\
\hline
3 & 31  & 4 & Theorem \ref{7be} (a) $+M(7)\le 4$\\
\hline
4 & 53 & 4 & Theorem \ref{7be} (b) $+M(7)\le 4$\\
\hline
5 & 47  & 4 & Theorem \ref{zeven} (b)\\
\hline
6 & 41  & 4 & Theorem\ref{zeven} (d)\\
\hline
\end{tabular}
\end{center}
Since $M(7;11)=M(7;17)=M(7;19)=4$ and $M(7;13)=3$ (the only cases not covered in Table 5), the proof
is completed. \qed\\

\noindent {\it Proof of Theorem} \ref{upto7}. Follows on combining Theorems \ref{driedrie}, \ref{5} and \ref{7}. 
\qed

\section{Computation of $M(11;q)$}
We have $M(11;q)\le M(11)=7$ (by Theorem \ref{ZZ} and Table 1). {}From \cite{GM} we recall the following result.
\begin{Thm}
\label{elf}
Let $q<r$ be primes such that $q\equiv 4({\rm mod~}11)$ and $r\equiv -3({\rm mod~}11)$. Let
$1\le \alpha \le q-1$ be the unique integer such that $11r\alpha\equiv 1({\rm mod~}q)$. Suppose
that
$q/33<\alpha\le (3q-1)/77$,
then $a_{11qr}(10+(6q-77\alpha)r)=-7$.
\end{Thm} 
\begin{Lem}
\label{Lem11_4}
Let $q$ be a prime such that $q\equiv 4({\rm mod~}11)$. For $q>37$, $M(11;q)=7$, and $M(11;37)=6$.
\end{Lem} 
{\it Proof}. By computation one finds that $M(11;37)=6$. Now assume $q>37$.
Notice that it is enough to show that $M(11;q)\ge 7$. For $q\ge 191$ the interval $I(q):=(q/33,(3q-1)/77]$ has length exceeding 1 and so contains at least one
integer $\alpha_1$. Then by the Chinese remainder theorem and Dirichlet's theorem we can find a prime
$r_1$ such that both $r_1\equiv -3({\rm mod~}11)$ and $11r_1\alpha_1\equiv 1({\rm mod~}q)$. 
Then we invoke Theorem \ref{elf} with $r=r_1$ and $\alpha=\alpha_1$.
It remains to
deal with the primes $59$ and $103$. One checks that both intervals $I(59)$ and $I(103)$ contain an integer and 
so we can proceed as in the case $q\ge 191$ to conclude the proof. \qed

\begin{Lem}
\label{Lem11_378}
Let $p=11$.\\
{\rm (a)} For $\ge 133$, $q\equiv 3 ({\rm mod~}11)$, $r^{-1}\equiv \frac{q-19}{2} ({\rm mod~}pq)$ and 
$k=q+7r\frac{(q-19)}{2}$ we have $a_{pqr}(k)=7$.\\
{\rm (b)} For $q\equiv 7 ({\rm mod~}11)$, $r^{-1}\equiv \frac{q+7}{2} ({\rm mod~}pq)$ and $k=6qr+4$ we have $a_{pqr}(k)=7$.\\
{\rm (c)} For $q\equiv 8 ({\rm mod~}11)$, $r^{-1}\equiv \frac{q-3}{2} ({\rm mod~}pq)$ and $k=6qr+4$ we have $a_{pqr}(k)=7$.
\end{Lem}

\noindent {\it Proof}. Apply Kaplan's lemma. \qed

\begin{Thm} 
\label{eleven}
For $q\ge 13$ we have
{\rm \begin{center}
\begin{tabular}{|c|c|c|c|c|c|c|c|c|c|c|}
\hline
$q({\rm mod~}11)$ & 1 & 2 & 3 & 4 & 5 & 6 & 7 & 8 & 9 & 10 \\
\hline
$M(11;q)$ & 6 & 6 & 7 & 7 & 6,7 & 6,7 & 7 & 7 & 6 & 6  \\
\hline
\end{tabular}
\end{center}}
\noindent except when $q\in \{17,23,37,43,47\}$. We have $M(11;17)=5$, $M(11;23)=3$, $M(11;37)=6$, $M(11; 43)=5$ and $M(11;47)=6$.
\end{Thm}
{\tt Remark 1}. If $q\equiv \pm 5({\rm mod~}11)$ and $q\ge 61$, then $M(p,q)\in \{6,7\}$.
We believe that $M(p;q)=6$.\\
{\tt Remark 2}. By Corollary \ref{minus1} and \ref{plus1} following Theorem \ref{newbound}, one infers that
$M(11;17)\le 5$, $M(11;23)\le 3$ and $M(11;43)\le 5$.\\

\noindent {\it Proof of Theorem} \ref{eleven}. We can most compactly prove this with a table.\\
\centerline{\bf TABLE 6}
\begin{center}
\begin{tabular}{|c|c|c|c|}
\hline
${\overline q}$ & $q_0$ & $M(11;q)$ & result\\
\hline
1 & 67 &  6 & Theorem \ref{zeven} (c)\\
\hline
2 & 13 & 6 & Theorem \ref{zeven} (a)\\
\hline
3 & 157 & 7 & Lemma \ref{Lem11_378} (a) $+M(11)\le 7$\\
\hline
4 & 59 & 7 & Lemma \ref{Lem11_4} \\
\hline
5 & 71  & 6,7 & Theorem \ref{7be} (a) $+M(11)\le 7$\\
\hline
6 & 61 & 6,7 & Theorem \ref{7be} (b) $+M(11)\le 7$\\
\hline
7 & 29 & 7 & Lemma \ref{Lem11_378} (b) $+M(11)\le 7$\\
\hline
8 & 19 & 7 & Lemma \ref{Lem11_378} (c) $+M(11)\le 7$\\
\hline
9 & 97  & 6 & Theorem \ref{zeven} (b)\\
\hline
10 & 109  & 6 & Theorem \ref{zeven} (d)\\
\hline
\end{tabular}
\end{center}
On directly computing the values of $M(p;q)$ not covered by the table, the proof is completed. \qed

\section{Computation for $p=19$}
By Theorem \ref{ZZ} we have $M(19)\le 2\cdot 19/3$ and hence $M(19)\le 12$. 
By Theorem \ref{oud} we find that $M(19;q)\ge 11$ for every $q\equiv 8({\rm mod~}19)$
and $q\ge 179$ and hence $M(19)\ge 11$.
Since $A(19\cdot 53\cdot 859)=12$, it follows that $M(19)=12$. The next result even shows that
$M(19;q)=M(19)$ for a positive fraction of the primes.
\begin{Thm}
\label{negentien}
We have $M(19)=12$.  Moreover, $M(19,q)=12$ if $q\equiv \pm 4({\rm mod~}19)$, with
$q>23$. Furthermore, $M(19;23)=11$.
\end{Thm}
The proof is an almost direct consequence of the following lemma.
\begin{Lem}
\label{19a}
Put $p=19$ and let $q\equiv \pm 4({\rm mod~}19)$ be a prime. Suppose there exists an
integer $a$ satysifying
\begin{equation}
\label{aaaa}
qa \equiv -1({\rm mod~}3)~{\rm and~}{q\over 6p}< a \le {5q-18\over 6p}.
\end{equation} 
Let $r>q$ be a prime satisfying $r(q-ap)\equiv 3({\rm mod~}pq)$. Then $a_{pqr}(7qr+q)=-12$, if
$q\equiv -4({\rm mod~}19)$, and $a_{19qr}(7qr+r)=-12$ if $q\equiv 4({\rm mod~}19)$.
\end{Lem}
{\it Proof}. Apply Kaplan's lemma. \qed\\

\noindent {\it Proof of Theorem} \ref{negentien}. For $q>90$ the interval in (\ref{aaaa}) is of length $>3$
and so contains an integer $a$ satisfying $qa \equiv -1({\rm mod~}3)$. It remains to deal with $q\in \{23,53,61\}$.
Computation shows that $M(19;23)=11$. For
$q=53$ and $q=61$ one finds an integer $a$ satisfying condition (\ref{aaaa}). \qed\\

\noindent {\it Proof of Theorem} \ref{exxpli}. By Theorem \ref{oud} and Dirichlet's theorem the claim follows for $p=13$.
Using Lemmas \ref{Lem11_4} and \ref{Lem11_378} the result follows for $p=11$.
On invoking Theorems \ref{upto7} and \ref{negentien}, the proof is then completed. \qed

\section{Small values of $M(p;q)$}
\label{close}
Typically if $M(p;q)$ is constant for all $q$ large enough with $q\equiv a({\rm mod~}d)$, then
$M(p;q)$ assumes a smaller value for some small $q$ in this progression. A (partial) explanation
of this phenomenon is provided in this section. We will show that if $ap+bq=1$ with $a$ and $b$
small in absolute value, then $M(p;q)$ is small. On the other hand we will show that $M(p;q)$ cannot
be truly small.
\begin{Prop}
\label{referee}
Let $2<p<q$ be odd primes. Then $M(p;q)\ge 2$.
\end{Prop}
\it Proof}. We say $\Phi_n(x)$ is flat if $A(n)=1$. ChunGang Ji \cite{Ji} proved that if
$p<q<r$ are odd prime and $2r\equiv \pm 1({\rm mod~}pq)$, then $\Phi_{pqr}(x)$ is flat iff
$p=3$ and $q\equiv 1({\rm mod~}3)$. It follows that $M(p;q)\ge 2$ for $p>3$. Now
invoke Theorem \ref{driedrie} to deal with the case $p=3$. \qed
\begin{Thm}
\label{newbound}
Let $2<p<q$ be odd primes and $\rho$ and $\sigma$ be the (unique) non-negative integers for which
$1+pq=(\rho+1)p+(\sigma+1)q$. Then
$$M(p;q)\le \begin{cases}p+\rho-\sigma & \mbox{ if }\rho\le \sigma; \\q+\sigma-\rho & \mbox{ if }\rho > \sigma.\end{cases}$$
\end{Thm}
\begin{cor}
\label{minus1}
Let $h,k$ be integers with $k>h$ and $q=(kp-1)/h$ a prime. If $p\ge k+h$, then $M(p;q)\le k+h$.
\end{cor}
\begin{cor}
\label{plus1}
Let $h,k$ be integers with $k>h$ and $q=(kp+1)/h$ a prime. If $p>h$ and $q> k+h$, then $M(p;q)\le k+h$.
\end{cor}
\noindent {\it Proof of Theorem }\ref{newbound}. Let us assume that $\rho\le \sigma$, the other case
being similar. Using Lemma \ref{flow} and Lemma \ref{binary} we infer that the number of $0\le m\le p-1$ with $b_{f(m)}=1$ is at most
$\rho+1$. Likewise the number of $m$ with $b_{f(m+q)}=-1$ is at most $p-1-\sigma$. By Kaplan's lemma it then
follows that $a_{pqr}(k)\le \rho+1+(p-1-\sigma)=p+\rho-\sigma$. Since the number of $0\le m\le p-1$ with
$b_{f(m)}=-1$ is at most $p-1-\sigma$ and the number of $m$ with $b_{f(m+q)}=1$ is at most $\rho+1$, we
infer that $a_{pqr}(k)\ge -(p+\rho-\sigma)$ and hence the result is proved. \qed

\begin{Thm}
\label{1modp}
Let $q\equiv 1({\rm mod~}p)$. Then
$$M(p;q)=\min\Big({q-1\over p}+1,{p+1\over 2}\Big).$$
\end{Thm}
{\it Proof}. For $p=3$ the result follows by Theorem \ref{driedrie}, so assume $p\ge 5$. Sister Beiter \cite{Beiter-1}, and
independently Bloom \cite{Bloom},  proved
that $M(p;q)\le (p+1)/2$ if $q\equiv \pm 1({\rm mod~}p)$ (alternatively we invoke Theorem \ref{upper1}).  By
Corollary \ref{plus1} we have $M(p;q)\le (q-1)/p+1$. By Lemma \ref{wilms2} the proof is then completed. \qed\\

\noindent Numerical experimentation suggests that in part (b) of
Theorem \ref{zeven} perhaps the condition $q>p^2/2$ can be dropped. By Theorem \ref{1modp} the condition
$q\ge (p-1)p/2+1$ in part (c) is optimal. In
part (d) we need $q\ge (p-1)p/2-1$, for otherwise $M(p;q)<(p+1)/2$ by Corollary \ref{minus1}.
\begin{Lem}
\label{driedrieX}
Let $p\ge 7$ be a prime such that $q=2p-1$ is also a prime. Let $r>q$ be a prime such
that $(p+q)r\equiv -2({\rm mod~}pq)$. Put $k=rq(p-1)/2+2p-pq$. Then $a_{pqr}(k)=3$.
\end{Lem}
{\it Proof}. Apply Kaplan's lemma. \qed\\

\noindent {\it Proof of Theorem} \ref{tweepeemineen}.  On combining Lemma \ref{driedrieX} with Corollary \ref{minus1}, one 
deduces that $M(p;2p-1)=3$ if $p\ge 5$ and $2p-1$ is a prime. \qed

\section{Conjectures, questions, problems}
The open problem that we think is the most interesting is Conjecture 1. Note that if one
could prove Conjecture 1 and getting an effective upper bound for the ternary
conductor $\mathfrak{f}_p$ (say $16p$) and an effective upper bound for the minimal ternary
prime (say $p^3$), then one has a finite procedure to compute $M(p)$. 
\begin{Prob}
Bachman {\rm \cite{B4}} introduced inclusion-exclusion polynomials. These polynomials generalize the ternary cyclotomic polynomials. Study $M(p;q)$ in this setting (here
$p$ and $q$ can be any coprime natural numbers), cf. Section 2 where we denoted this function by
$M'(p;q)$. For example, using \cite[Theorem 3]{B4} by an argument similar to
that given in Proposition \ref{finitep} it is easily seen that there is a finite procedure
to compute $M'(p;q)$.
\end{Prob}
\begin{Prob}
The analogue of $M(p;q)$ for inverse cyclotomic polynomials, see {\rm \cite{Moree}}, can 
be defined. Study it.
\end{Prob}
\begin{Qu}
Can one compute the average value of $M(p;q)$, that is does the limit
$$\lim_{x\rightarrow \infty}{1\over \pi(x)}\sum_{p<q\le x}M(p;q)$$
exist and if yes, what is its value?
\end{Qu}
\begin{Qu}
Is Theorem {\rm \ref{exxpli}} still true if we put $\delta(13)=1/3$ and cross out the words `a subset having'?
\end{Qu}
\begin{Qu}
If $q>p$ is prime and $q\equiv -2({\rm mod~}p)$, then do we have $M(p;q)=(p+1)/2$?
\end{Qu}
\begin{Qu}
Suppose that $p>11$ is a prime.\\
If $6p - 1$ is prime, then do we have $M(p, 6p - 1) = 7$?\\
If $(5p - 1)/2$ is prime, then do we have $M(p, (5p - 1)/2) = 7$?\\
If $(5p + 1)/2$ is prime then do we have $M(p, (5p + 1)/2) = 7$?\\
Find more similar results.
\end{Qu}
\begin{Qu}
Given an integer $k\ge 1$, does there exist $p_0(k)$ and a function $q_k(p)$ such that if
$q\equiv 2/(2k+1)({\rm mod~}p)$, $q\ge q_k(p)$ and
$p\ge p_0(k)$, then $M(p;q)=(p+2k+1)/2$?
\end{Qu}
\begin{Qu}
Is it true that $M(11;q)=6$ for all large enough $q$ satisfying $q\equiv \pm 5({\rm mod~}11)$ ?
If so one can finish the computation of $M(11;q)$.
\end{Qu}
\begin{Qu}
\label{zeven77}
Is it true that for $q$ sufficiently large the values of $M(13;q)$, $M(17;q)$, $M(19;q)$ and $M(23;q)$ are given by the following tables?
\end{Qu}
{\rm
\begin{center}
\begin{tabular}{|c|c|c|c|c|c|c|c|c|c|c|c|c|}
\hline
$q({\rm mod~}13)$ & 1 & 2 & 3 & 4 & 5 & 6 & 7 & 8 & 9 & 10 & 11 & 12\\
\hline
$M(13;q)$ & 7 & 7 & 7 & 8 & 8 & 7 & 7 & 8 & 8 & 7 & 7 & 7 \\
\hline
\end{tabular}
\end{center}

\begin{center}
\begin{tabular}{|c|c|c|c|c|c|c|c|c|c|c|c|c|c|c|c|c|}
\hline
$q({\rm mod~}17)$ & 1 & 2 & 3 & 4 & 5 & 6 & 7 & 8 & 9 & 10 & 11 & 12 & 13 & 14 & 15 & 16\\
\hline
$M(17;q)$ & 9 & 9 & 9 & 10 & 10 & 9 & 10 & 9 & 9 & 10 & 9 & 10 & 10 & 9 & 9 & 9\\
\hline
\end{tabular}
\end{center}

\begin{center}
\begin{tabular}{|c|c|c|c|c|c|c|c|c|c|c|c|c|c|c|c|c|c|c|}
\hline
$q({\rm mod~}19)$ & 1 & 2 & 3 & 4 & 5 & 6 & 7 & 8 & 9\\
\hline
$M(19;q)$ & 10 & 10 & 10 & 12 & 11 & 9 & 11 & 11 & 10\\
\hline
\hline
$q({\rm mod~}19)$ & 10 & 11 & 12 & 13 & 14 & 15 & 16 & 17 & 18\\
\hline
$M(19;q)$ & 10 & 11 & 11 & 9 & 11 & 12 & 10 & 10 & 10\\
\hline
\end{tabular}
\end{center}

\begin{center}
\begin{tabular}{|c|c|c|c|c|c|c|c|c|c|c|c|c|c|c|c|c|c|c|c|c|}
\hline
$q({\rm mod~}23)$ & 1 & 2 & 3 & 4 & 5 & 6 & 7 & 8 & 9 & 10 & 11\\
\hline
$M(23;q)$ & 12 & 12 & 12 & 14 & 14 & 11 & 13 & 11 & 14 & 13 & 12\\
\hline
\hline
$q({\rm mod~}23)$ & 12 & 13 & 14 & 15 & 16 & 17 & 18 & 19 & 20 & 21 & 22\\
\hline
$M(23;q)$ & 12 & 13 & 14 & 11 & 13 & 11 & 14 & 14 & 12 & 12 & 12\\
\hline
\end{tabular}
\end{center}}
The next question is raised by the referee of this paper.
\begin{Qu}
Suppose that for all sufficiently large primes 
$q\equiv q_0({\rm mod~}{\mathfrak f}_p)$ we have
$M(p;q)<M(p)$. Is it possible to prove that $M(p;q)<M(p)$ for every prime
$q\equiv q_0({\rm mod~}{\mathfrak f}_p)$?
\end{Qu}
\begin{Qu}
For a given prime $p$, let $m(p)$ denote $\lim \inf M(p;q)$, with $q>p$.
Determine $m(p)$. Is it true that $\lim_{p\rightarrow \infty}m(p)/p=c$ for some constant $c>0$?
\end{Qu}
By Proposition \ref{referee} we have $m(p)\ge 2$ for $p>2$.
Note that the results in this paper imply that $m(p)=(p+1)/2$ for $2<p\le 11$. If the answer to
Question \ref{zeven77} is yes, then $m(p)=(p+1)/2$ for $2<p\le 17$ and $m(p)=(p-1)/2$ for
$19\le p \le 23$.\\ 
(The issue of lower bounds for $M(p;q)$ was raised by the referee.)\\

\noindent {\tt Acknowledgement}. The second author likes to thank the MPIM interns (in June
2008) N. Baghina, C. Budde, B. J\"uttner and D. Sullivan for their (computer) assistance in computing
tables used in the proofs of Lemma \ref{19a} ($p=19$). For these tables see \cite{GMW}.
However, the bulk of the paper was written whilst the third author was during two months in 2010 an 
intern at MPIM under the guidance of the second author. The third author would like to thank the MPIM for the 
possibility to do an internship and for the nice research atmosphere. He also thanks the second author for his 
mentoring and for having a sympathetic ear for any questions.\\
\indent Finally thanks are due to G. Bachman for some helpful remarks, and the referee who spent quite a bit
of time writing a very extensive referee report, which led to many improvements over the original submission.

\medskip\noindent {\footnotesize 12 bis rue Perrey, 31400 Toulouse, France.\\
e-mail: {\tt galloty@orange.fr}}\\

\medskip\noindent {\footnotesize Max-Planck-Institut f\"ur Mathematik,\\
Vivatsgasse 7, D-53111 Bonn, Germany.\\
e-mail: {\tt moree@mpim-bonn.mpg.de}\\

\medskip\noindent {\footnotesize Sterbeckerstrasse 21, 58579 Schalksm\"uhle, Germany\\
e-mail: {\tt Robert.wilms@rub.de}}

\vskip 5mm


{\small
\begin{thebibliography}{99}
\bibitem{B1} G. Bachman, On the coefficients of ternary cyclotomic polynomials, 
{\it J. Number Theory} {\bf 100} (2003), 104--116.
\bibitem{B2} G. Bachman, Ternary cyclotomic polynomials with an optimally large set of 
coefficients, {\it Proc. Amer. Math. Soc.}  {\bf 132} (2004), 1943--1950. 
\bibitem{B3} G. Bachman, Flat cyclotomic polynomials of order three, 
{\it Bull. London Math. Soc.} {\bf 38} (2006), 53--60.
\bibitem{B4} G. Bachman, On ternary inclusion-exclusion polynomials, {\it Integers} {\bf 10} (2010), 
A48, 623--638. 
\bibitem{BM} G. Bachman and P. Moree, On a class of ternary
inclusion-exclusion polynomials, {\it Integers} {\bf 11} (2011), A8, 14 pp.
\bibitem{Bang} A.S. Bang, Om Ligningen $\varphi_n(x)=0$, {\it Nyt Tidsskrift for Mathematik (B)}
{\bf 6} (1895), 6--12. 
\bibitem{Beiter-1} M. Beiter, Magnitude of the coefficients of the cyclotomic 
polynomial $F\sb{pqr}\,(x)$, {\it Amer. Math. Monthly}  {\bf 75} (1968), 370--372.
\bibitem{Beiter-2} M. Beiter, Magnitude of the coefficients of the cyclotomic 
polynomial $F\sb{pqr}$. II, {\it Duke Math. J.} 
{\bf 38} (1971), 591--594.
\bibitem{Beiter-3} M. Beiter, Coefficients of the cyclotomic polynomial 
$F_{3qr}(x)$, {\it Fibonacci Quart.}  {\bf 16}  (1978), 302--306.
\bibitem{Bloom} D.M. Bloom, On the coefficients of the cyclotomic polynomials, 
{\it Amer. Math. Monthly} {\bf 75} (1968), 372--377. 
\bibitem{BZ} B. Bzd{\c e}ga, Bounds on ternary cyclotomic coefficients, 
 {\it Acta Arith.}  {\bf 144}  (2010), 5–-16. 
\bibitem{CGM} C. Cobeli, Y. Gallot and P. Moree, unpublished manuscript.
\bibitem{DFI} W. Duke, J.B. Friedlander and H. Iwaniec, Equidistribution of roots of a quadratic 
congruence to prime moduli, {\it Ann. of Math.} (2)  {\bf 141}  (1995),  423--441. 
\bibitem{jessica} J. Fintzen, Cyclotomic polynomial coefficients $a(n,k$) with $n$ and $k$ in 
prescribed residue classes, {\it J. Number Theory} {\bf 131} (2011), 1852--1863. 	
\bibitem{GM} Y. Gallot and P. Moree,  Ternary cyclotomic polynomials having a large
coefficient, {\it  J. Reine Angew. Math.}  {\bf 632}  (2009), 105--125.
\bibitem{GM2} Y. Gallot and P. Moree, Neighboring ternary cyclotomic coefficients differ by at 
most one, {\it J. Ramanujan Math. Soc.} {\bf 24} (2009), 235--248.
\bibitem{GMW} Y. Gallot, P. Moree and R. Wilms, The family of ternary cyclotomic polynomials with one free prime,
MPIM-preprint 2010-11, pp. 32.
\bibitem{Ji} C. Ji, A specific family of cyclotomic polynomials of order three, 
{\it Sci. China Math.} {\bf 53} (2010), 2269--2274.
\bibitem{Kaplan} N. Kaplan, Flat cyclotomic polynomials of order three, {\it J. Number Theory} 
{\bf 127} (2007), 118--126.
\bibitem{LL} T.Y. Lam and K.H. Leung, On the cyclotomic polynomial $\Phi_{pq}(X)$, 
{\it Amer. Math. Monthly} {\bf 103} (1996), 562--564. 
\bibitem{Emma} E. Lehmer, On the magnitude of the coefficients of the cyclotomic 
polynomials, {\it Bull. Amer. Math. Soc.} {\bf 42} (1936), 389--392. 
\bibitem{HM} H. M\"oller, \"Uber die Koeffizienten des $n$-ten Kreisteilungspolynoms, 
{\it Math. Z.} {\bf 119} (1971), 33--40. 
\bibitem{Moree} P. Moree, Inverse cyclotomic 
polynomials, {\it J. Number Theory} {\bf 129} (2009), 667--680.
\bibitem{Thanga} R. Thangadurai, On the coefficients of 
cyclotomic polynomials, {\it Cyclotomic fields and related topics} 
(Pune, 1999), 311--322, Bhaskaracharya
  Pratishthana, Pune, 2000. 
\bibitem{ZZ1} J. Zhao and X. Zhang, Coefficients of ternary cyclotomic polynomials, 
{\it J. Number Theory} {\bf 130} (2010), 2223-2237.
\bibitem{ZZ2} J. Zhao and X. Zhang, A proof of the Corrected Beiter conjecture, 
arXiv:0910.2770.
\end{thebibliography}}
\end{document}